\documentclass{amsart}
\usepackage{latexsym}
\usepackage{amssymb}

\usepackage[pctex32]{graphics}
\newcommand{\real}{\mathbb R}

\begin{document}
\title[VECTOR VALUED FUNCTIONS NOT CONSTANT ON CONNECTED SETS]{VECTOR VALUED FUNCTIONS NOT CONSTANT ON CONNECTED SETS OF CRITICAL POINTS}
\author{Azat Ainouline}
\address{}
\email{azat.ainouline@3web.net}
\begin{abstract}
Whitney type examples of maps $f\in C^k(\real^m,\real^n)$ for a maximal possible real $k$, and multidimensional space-filling curves with special properties are constructed.
\end{abstract}
\maketitle

\newtheorem{theorem}{Theorem}
\newtheorem{lemma}{Lemma}[section]
\newcommand{\diam}{{\rm diam}}
\newtheorem{definition}{Definition}[section]
\numberwithin{equation}{section}

\section{Introduction.}
In 1935 Hassler Whitney \cite{Whitney} published his example of a $C^1$ function $f:\real^2\rightarrow\real$~~not constant on a connected set of critical points.
The following theorem is the main result of this paper.

\begin{theorem}\label{theorem 1}\rm
For any  $n,m\in \mathbb N$ there exist a map $p:[0,1]^m\rightarrow[0,1]^n$, contained in $C^k$ for all real $k<\frac{m}{n}$, and a connected set $E\subseteq[0,1]^m$, such that every partial derivative of $p$ of order $<\frac{m}{n}$ vanishes on $E$ and $p(E) =[0,1]^n.$
\end{theorem}

\noindent
Easy corollary from the Theorem \ref{theorem 1} is:
\begin{theorem}\label{theorem 2}\rm
Let $n,m,r$ be non-negative integer numbers,~ $m>n>r$,~then there exists a map $p:\real^m\rightarrow \real^n$, contained in $C^k$ for all real $k<\frac{m-r}{n-r}$, and a connected subset $E$  of points of rank $r$ such that $p(E)$ contains an open set.
\end{theorem}

\noindent
Similar to the Theorem \ref{theorem 2} results, but without connectedness of $E$, were obtained by R.Kaufman \cite{Kaufman}, S.M.Bates \cite{Bates 1},\cite{Bates 2}, A.Norton \cite{Norton}.\\ 

\noindent
The following theorem, being a weak version of the Bates' theorem  \cite{Bates 2}, shows the sharpness of the $k$ in Theorems \ref{theorem 1},\ref{theorem 2}:\\
{\bf Theorem (Bates)}
Let $n,m,r$ be nonnegative integers satisfying $m>n>r$, and define $s=(m-r)/(n-r)$.~ If $E$ is a set of rank $r$ for $f:\real^m\rightarrow\real^n$~ and $f\in C^s$,
then $f(E)$~ has Lebesgue measure zero in $\real^n$.\\

Whitney indicated in same paper \cite{Whitney} how to generalize his construction for higher dimensions (ie. $f\in C^k(\real^m,\real^1)$~ for~$m,k\in\mathbb N,~~k<m$). He wrote for the case $m=3,k=2$~~"Let $Q$ be a cube of side 1. Let $Q_0,\dots,Q_7$~ be cubes of side 2/5 arranged in $Q$ so that $Q_i$ is adjacent to $Q_{i-1}$". But practically this "adjacency" is not such easy to make. The problem deeply relates to the problem of construction of multidimensional space-filling curves with special properties, which is quite complicated itself. One can find the properties that fit Whitney's example for cases $m=3$ in Sagan's curve $f:[0,1]\rightarrow[0,1]^3$ ~(see \cite{Sagan 2}) published in 1993, and for case $m=4$ in Steinhause' curve $f:[0,1]\rightarrow[0,1]^4$ ~ (see \cite{Sagan}). In Theorems \ref{big lemma},\ref{corollary} of this paper the author constructs the space-filling function $f:[0,1]\rightarrow[0,1]^m$ for any $m\in\mathbb N$ with the special "adjacency", existence of which Whitney assumed. \\   

 Such a type of space-filling curves appears to be a very useful tool in solution of some problems involving Euclidian space and $C^k$ functions. One of such applications is the Whitney's type examples constructed in this paper. 

\subsubsection{Preliminary results.}

\begin{definition}
\rm If $f:\real^m\rightarrow\real^n$, and $\lambda\in(0,1]$, we define $\lambda$-partial derivatives $f_1^{(\lambda)},\dots,f_m^{(\lambda)}$~ by the formula:
\begin{align*}
f_i^{(\lambda)}(a) =\lim\limits_{t\to 0}sign(t)\frac{f(a_1,\dots,a_{i-1},a_i+t,a_{i+1},\dots,,a_n)-f(a)}{|t|^\lambda}
\end{align*}
for $a=(a_1,\dots,a_m)\in \real^m$.  If all $\lambda$-partial derivatives are continuous on some $M\subseteq\real^m$ 
 we say that $f\in C^\lambda $ on $M$.
\end{definition}

\begin{definition}
\rm For $k\in \real^+$~~a function~~$ f:M\subseteq\real^m\rightarrow \real^n$~~is a $C^k$-function (or $f\in C^k$),~~if $f\in C^{[k]}$ and every $[k]^{th}$ partial derivative of $f$ on $M$ is $C^{k-[k]}$-function,~where $[k]$~is the integer part of $k$. If $f\in C^k$~ for every ~$k<k_0$,~ we will write $f\in C^{<k_0}.$
\end{definition}

We begin by setting ~$K^n_0=\{Q^n_1\}$,~where $Q^n_1=[0,1]^n$ is the closed cube in $\mathbb R^n$ with side length 1.
\newline 
In general, having constructed the cubes of $K^n_{s-1}$,  divide each $Q^n_{1,i_1,i_2,....,i_{s-1}}\in K^n_{s-1}$ into $2^n$ closed cubes of side $\frac{1}{2^s}$,  and let $K^{n}_{s}$ be the set of all these cubes. More precisely we will write 

$$K^{n}_{s}=\{ Q^n_{1,i_1,i_2,...,i_{s-1},i_s}~~;~~ Q^n_{1,i_1,i_2,...,i_{s-1},i_s} \subseteq Q^n_{1,i_1,i_2,...i_{s-1}}\in K^n_{s-1}, 1\leqslant i_s\leqslant 2^n \}.$$

\noindent We also define:
\begin{itemize}
\item
$K^n=\bigcup \limits_{s\in N} K^n_s$ ~~(note: $K^n$ is defined for $\mathbb R^n$);
\item     
$S(Q^n_{1,i_1,i_2,...,i_s})$ = $\frac{1}{2^s}$ ~- the length of a side of $Q^n_{1,i_1,i_2,...,i_s}$;
\item
$meas(\delta)$ be the Lebesgue measure of $\delta\subseteq\mathbb R^n$.
\end{itemize}

\noindent The main goal of this paragraph is to construct for any $n\in\mathbb N$ a continuous space filling curve $$f_n :~[0,1] \stackrel{onto}{\longrightarrow} [0,1]^n$$ with  special properties:

\begin{align}
\label{cubes preserving 1}
&\text{if}~\alpha\subseteq [0,1]~~~\text{and~for~some}~~s\in\mathbb N~~~&\alpha \in K^1_{n\cdot s}~&~ \text{then}~~f_n(\alpha) \subseteq \delta~~\text{for~some}~\delta\in K^n_s\\
&\text{if}~\delta\subseteq [0,1]^n~~\text{and~for~some}~~s\in\mathbb N~~~&\delta \in 
K^n_s~&~ \text{then}~~f_n^{-1}(int(\delta)) \subseteq\alpha~~\text{for~some}~\alpha\in K^1_{n\cdot s}
\label{cubes preserving 2}
\end{align}

\noindent
where~~$int(\delta)$~~is~the~set~of~interior~points~of~~$\delta$.

\begin{definition}
\rm We will call a function $f_n:~[0,1]\rightarrow[0,1]^n$ with the properties \eqref{cubes preserving 1},\eqref{cubes preserving 2} {\bf cubes preserving}.
\end{definition}
Note that a continuous cubes preserving function $f_n$  is a space-filling and measure preserving function with  a property:
$$\text{if}~\alpha\subseteq[0,1]~\text{and~for~some~} s\in\mathbb N~~\alpha\in K^1_{n\cdot s}~~\text{then}~~f_n(\alpha)=\delta~\text{for~some}~\delta\in K^n_s.$$ 

\begin{definition}
\rm Extending Sagan's definition of measure preserving function \cite{Sagan} to dimension $n$, we will call a function $q:[0,1]\rightarrow [0,1]^n$ {\bf measure preserving}, if for every $P\subseteq [0,1]^n~~ meas(q^{-1}(P)) = meas(P)$.
\end{definition}

Once we have a cube preserving function $q:[0,1]\rightarrow[0,1]^n$ for some $n\in\mathbb N$, it allows us to build a linear order on $K^n_s$ with respect to the $q$  for $s\in\mathbb N$:
\begin{align}\label{linear order}
&\forall \delta,\delta^\prime\in K^n_s~~~~\delta\prec_q\delta^\prime~~~~\text{if}~~~~\forall x\in q^{-1}(int(\delta)),x^\prime\in q^{-1}(int(\delta^\prime))~~~~x<x^\prime.
\end{align}
So that we can enumerate the elements of $K^n_s$ as $\{Q^n_{1,i_1,i_2,\dots,i_j,\dots,i_s};~1\leqslant i_j\leqslant 2^n\}$, where $Q^n_{1,i_1,i_2,\dots,i_s}~\prec_q~Q^n_{1,i_1^\prime,i_2^\prime,\dots,i_s^\prime}$ if $\nexists j_0~(1\leqslant j_0\leqslant s)$~ such that $\forall j<j_0~~i_j\geqslant i_j^\prime$ and $i_{j_0}^\prime<i_{j_0}.$ Also note that $Q^n_{1,i_1,i_2,\dots,i_{s-1},i_s}\subset Q^n_{1,i_1,i_2,\dots,i_{s-1}}~\forall i_s~~~(1\leqslant i_s\leqslant 2^n).$
Thereby we have enumerated all elements of the set $K^n$. Let us designate the enumeration by $\langle K^n,\prec_q\rangle$.\\

Note that cubes preserving functions were already constructed for $n\leqslant 4$. For example:  
\begin{itemize}
\item
In case $n=2$: the Hilbert curve $f_H :~[0,1]\rightarrow [0,1]^2$ and the Peano curve $f_P :~[0,1]\rightarrow [0,1]^2$. 
\item
In case $n=3$: the Sagan's generalization of the Hilbert curve $f_S:~[0,1]\rightarrow [0,1]^3$. 
\item
In case $n=4$: the Steinhause Space-Filling curve $f_{SH}:~[0,1]\rightarrow [0,1]^4$ with its coordinate functions: $\psi\psi,\psi\varphi,\varphi\psi,\varphi\varphi$,~~ where $\psi,\varphi$ - are the coordinate functions of the Hilbert curve.  
\end{itemize}

\noindent
We now will prove that one can construct a cubes preserving function for any $n\in\mathbb N$.

\begin{lemma}
\label{small lemma}
\rm
Let $E_1,E_2$ be copies of $\mathbb R$. 
\newline $\forall n\in\mathbb N~~~\text{there exists}~~~ \text{continuous}~S_n:~[0,1]\stackrel{\rm onto}{\longrightarrow}[0,1]^2\subseteq E_1\times E_2$  such that 
\begin{enumerate}
\item \label{Sn property 1}
 If $\alpha\in K^1_{n\cdot s}$, then $S_n(\alpha)\subseteq  \alpha^\prime\times\alpha^{\prime\prime}$,  where $\alpha^\prime\in K^1_{(n-1)\cdot  s},~~\alpha^{\prime\prime}\in K^1_s,~\alpha^\prime\subset E_1,~\alpha^{\prime\prime}\subset E_2$
\item \label{Sn property 2}
 If $\alpha^\prime\times\alpha^{\prime\prime}\subseteq[0,1]^2$  such that $\alpha^\prime\subseteq E_1,~~\alpha^{\prime\prime}\subseteq E_2,~\alpha^\prime\in K^1_{(n-1)\cdot s},~\alpha^{\prime\prime}\in K^1_s$,~~then $S^{-1}_n(int(\alpha^\prime\times\alpha^{\prime\prime}))\subseteq\alpha\in K^1_{n\cdot s}$.
\end{enumerate}
\end{lemma}

\noindent
{\bf Proof.}
Let us define for every $n \geqslant 2$ a function $S_n$ as follows:

If the interval $[0,1]$ can be mapped continuously onto the square $[0,1]^2$, then after partitioning $[0,1]$ into $2^n$ congruent subintervals and $[0,1]^2$ into $2^n$ congruent subrectangles with sides $\frac{1}{2^{n-1}},~\frac{1}{2^1}$,  each subinterval can be mapped continuously onto one of the subrectangles. 

Next, each subinterval is, in turn, partitioned into $2^n$ congruent subintervals, and each subrectangles into $2^n$ congruent subrectangles with sides $\frac{1}{2^{(n-1)\cdot 2}}, ~\frac{1}{2^2}$ and the argument is repeated. If this is carried on indefinitely, $[0,1]$ and $[0,1]^2$ are partitioned into $2^{ns}$ congruent replicas each with sides $\frac{1}{2^{(n-1)s}},~ \frac{1}{2^s}$ for $s\in\mathbb N$. Fig.\ref{blank form} shows how we divide each rectangle(square in first iteration) into subrectangles. 

\begin{figure}[ht]
\begin{minipage}[t]{12cm}
\setlength{\unitlength}{0.08cm}
\begin{picture}(40,50)
\put(50,0)   {\framebox(40,50)}
\put(50,0)   {\framebox(20,10){1}}
\put(70,0)   {\framebox(20,10){2}}
\put(50,10)  {\framebox(20,10){4}}
\put(70,10)  {\framebox(20,10){3}}
\put(50,20)  {\framebox(20,10){5}}
\put(70,20)  {\framebox(20,10){6}}
\put(50,30)  {\framebox(20,10){$\dots$}}
\put(70,30)  {\framebox(20,10){$\dots$}}
\put(50,40)  {\framebox(20,10){$2^n$}}
\put(70,40)  {\framebox(20,10){$2^{n-1}$}}
\end{picture}
\caption{Each iteration $s$ divides subrectangle into $2^n$ subsubrectangles with sides $\frac{1}{2^{(n-1)\cdot s}},~\frac{1}{2^s}$}
\label{blank form}
\end{minipage}
\hfill
\end{figure}

We need to demonstrate that the subrectangles can be arranged so that adjacent subintervals correspond to adjacent subrectangles with an edge in common, and so that the inclusion relationships are presented, i.e. if a rectangle corresponds to an interval, then its subrectangles correspond to the subintervals of that interval. 

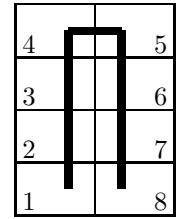
\begin{figure}[htb]
\begin{minipage}[t]{6cm}
\setlength{\unitlength}{0.07cm}
\begin{picture}(30,40)
\thinlines
\put(0,0)     {\framebox(30,40)}
\put(0,10)    {\line(1,0){30}}
\put(0,20)    {\line(1,0){30}}
\put(0,30)    {\line(1,0){30}}
\put(15,0)    {\line(0,1){40}}

\put(40,0)    {\framebox(30,40)}
\put(40,10)   {\line(1,0){30}}
\put(40,20)   {\line(1,0){30}}
\put(40,30)   {\line(1,0){30}}
\put(55,0)    {\line(0,1){40}}

\put(0,0)    {\makebox (5,5){1}}
\put(25,0)   {\makebox (5,5){2}}
\put(0,10)   {\makebox (5,5){4}}
\put(25,10)  {\makebox (5,5){3}}
\put(0,20)   {\makebox (5,5){5}}
\put(25,20)  {\makebox (5,5){6}}
\put(0,30)   {\makebox (5,5){8}}
\put(25,30)  {\makebox (5,5){7}}
\put(40,0)   {\makebox (5,5){2}}
\put(65,0)   {\makebox (5,5){1}}
\put(40,10)  {\makebox (5,5){3}}
\put(65,10)  {\makebox (5,5){4}}
\put(40,20)  {\makebox (5,5){6}}
\put(65,20)  {\makebox (5,5){5}}
\put(40,30)  {\makebox (5,5){7}}
\put(65,30)  {\makebox (5,5){8}}
\linethickness{1mm}
\put(10,5)   {\line(1,0){10}}
\put(20,5)   {\line(0,1){10}}
\put(10,15)  {\line(1,0){10}}
\put(10,15)  {\line(0,1){10}}
\put(10,25)  {\line(1,0){10}}
\put(20,25)  {\line(0,1){10}}
\put(10,35)  {\line(1,0){10}}

\put(50,5)   {\line(1,0){10}}
\put(50,5)   {\line(0,1){10}}
\put(50,15)  {\line(1,0){10}}
\put(60,15)  {\line(0,1){10}}
\put(50,25)  {\line(1,0){10}}
\put(50,25)  {\line(0,1){10}}
\put(50,35)  {\line(1,0){10}}
\end{picture}
\caption{Peano Method ("P")}
\label{Peano Method}
\end{minipage}
\hfill
\begin{minipage}[t]{6cm}
\setlength{\unitlength}{0.07cm}
\begin{picture}(30,40)
\put(10,0)   {\framebox(30,40)}
\put(10,10)    {\line(1,0){30}}
\put(10,20)    {\line(1,0){30}}
\put(10,30)    {\line(1,0){30}}
\put(25,0)    {\line(0,1){40}}

\put(50,0)    {\framebox(30,40)}
\put(50,10)   {\line(1,0){30}}
\put(50,20)   {\line(1,0){30}}
\put(50,30)   {\line(1,0){30}}
\put(65,0)    {\line(0,1){40}}

\put(10,0)   {\makebox (5,5){4}}
\put(35,0)   {\makebox (5,5){5}}
\put(10,10)  {\makebox (5,5){3}}
\put(35,10)  {\makebox (5,5){6}}
\put(10,20)  {\makebox (5,5){2}}
\put(35,20)  {\makebox (5,5){7}}
\put(10,30)  {\makebox (5,5){1}}
\put(35,30)  {\makebox (5,5){8}}
\put(50,0)   {\makebox (5,5){1}}
\put(75,0)   {\makebox (5,5){8}}
\put(50,10)  {\makebox (5,5){2}}
\put(75,10)  {\makebox (5,5){7}}
\put(50,20)  {\makebox (5,5){3}}
\put(75,20)  {\makebox (5,5){6}}
\put(50,30)  {\makebox (5,5){4}}
\put(75,30)  {\makebox (5,5){5}}
\linethickness{1mm}
\put(20,5)   {\line(1,0){10}}
\put(20,5)   {\line(0,1){30}}
\put(30,5)   {\line(0,1){30}}

\put(60,35)  {\line(1,0){10}}
\put(60,5)   {\line(0,1){30}}
\put(70,5)   {\line(0,1){30}}

\end{picture}
\caption{Hilbert Method ("H")}
\label{Hilbert Method}
\end{minipage}
\hfill
\end{figure}

We will use here a combination of two different methods to construct these space-filling curves:
\begin{enumerate}
\item
 The first method is based on idea of Peano \cite{Sagan}.  Fig. \ref{Peano Method} shows the way of construction the curve by Peano method. For the future use, we designate  this method as "P".
\item
 The second method is based on idea of Hilbert \cite{Sagan}. Fig. \ref{Hilbert Method} shows the ways of construction the curve by Hilbert method. We designate this method as "H".  
\end{enumerate}
Note: For both methods, the choice of which of two ways will be used depends only on the choice of the start point, and  disposition of the end point follows the parity of $2^{n-1}$ and does not depend on $n$.

\noindent
To create the next iteration curve, we will give the means of how to present each subectangle from the previous iteration (see Fig.\ref{P-continue},\ref{H-continue}).

\begin{figure}[htb]
\begin{minipage}[t]{6cm}
\setlength{\unitlength}{0.07cm}
\begin{picture}(30,40)
\thinlines
\put(0,0)    {\framebox(30,40)}
\put(0,10)    {\line(1,0){30}}
\put(0,20)    {\line(1,0){30}}
\put(0,30)    {\line(1,0){30}}
\put(15,0)    {\line(0,1){40}}

\put(40,0)    {\framebox(30,40)}
\put(40,10)    {\line(1,0){30}}
\put(40,20)    {\line(1,0){30}}
\put(40,30)    {\line(1,0){30}}
\put(55,0)     {\line(0,1){40}}

\put(0,0)    {\makebox (5,5){H}}
\put(25,0)   {\makebox (5,5){P}}
\put(0,10)   {\makebox (5,5){H}}
\put(25,10)  {\makebox (5,5){P}}
\put(0,20)   {\makebox (5,5){H}}
\put(25,20)  {\makebox (5,5){P}}
\put(0,30)   {\makebox (5,5){H}}
\put(25,30)  {\makebox (5,5){P}}
\put(40,0)   {\makebox (5,5){P}}
\put(65,0)   {\makebox (5,5){H}}
\put(40,10)  {\makebox (5,5){P}}
\put(65,10)  {\makebox (5,5){H}}
\put(40,20)  {\makebox (5,5){P}}
\put(65,20)  {\makebox (5,5){H}}
\put(40,30)  {\makebox (5,5){P}}
\put(65,30)  {\makebox (5,5){H}}
\linethickness{1mm}
\put(10,5)   {\line(1,0){10}}
\put(20,5)   {\line(0,1){10}}
\put(10,15)  {\line(1,0){10}}
\put(10,15)  {\line(0,1){10}}
\put(10,25)  {\line(1,0){10}}
\put(20,25)  {\line(0,1){10}}
\put(10,35)  {\line(1,0){10}}

\put(50,5)   {\line(1,0){10}}
\put(50,5)   {\line(0,1){10}}
\put(50,15)  {\line(1,0){10}}
\put(60,15)  {\line(0,1){10}}
\put(50,25)  {\line(1,0){10}}
\put(50,25)  {\line(0,1){10}}
\put(50,35)  {\line(1,0){10}}
\end{picture}
\caption{The means of moving to the next iteration in method "P"}
\label{P-continue}
\end{minipage}
\hfill
\begin{minipage}[t]{6cm}
\setlength{\unitlength}{0.07cm}
\begin{picture}(30,40)
\put(10,0)   {\framebox(30,40)}
\put(10,10)    {\line(1,0){30}}
\put(10,20)    {\line(1,0){30}}
\put(10,30)    {\line(1,0){30}}
\put(25,0)     {\line(0,1){40}}

\put(50,0)    {\framebox(30,40)}
\put(50,10)    {\line(1,0){30}}
\put(50,20)    {\line(1,0){30}}
\put(50,30)    {\line(1,0){30}}
\put(65,0)     {\line(0,1){40}}

\put(10,0)   {\makebox (5,5){H}}
\put(35,0)   {\makebox (5,5){H}}
\put(10,10)  {\makebox (5,5){P}}
\put(35,10)  {\makebox (5,5){P}}
\put(10,20)  {\makebox (5,5){P}}
\put(35,20)  {\makebox (5,5){P}}
\put(10,30)  {\makebox (5,5){P}}
\put(35,30)  {\makebox (5,5){P}}
\put(50,0)   {\makebox (5,5){P}}
\put(75,0)   {\makebox (5,5){P}}
\put(50,10)  {\makebox (5,5){P}}
\put(75,10)  {\makebox (5,5){P}}
\put(50,20)  {\makebox (5,5){P}}
\put(75,20)  {\makebox (5,5){P}}
\put(50,30)  {\makebox (5,5){H}}
\put(75,30)  {\makebox (5,5){H}}
\linethickness{1mm}
\put(20,5)   {\line(1,0){10}}
\put(20,5)   {\line(0,1){30}}
\put(30,5)   {\line(0,1){30}}

\put(60,35)  {\line(1,0){10}}
\put(60,5)   {\line(0,1){30}}
\put(70,5)   {\line(0,1){30}}
\end{picture}
\caption{The means of moving to the next iteration in method "H"}
\label{H-continue}
\end{minipage}
\hfill
\end{figure}

\noindent
And finally in Fig.\ref{P-second iteration},\ref{H-second iteration} we indicate how this process is to be carried out for the next iteration. 
\newline
\begin{figure}[htb]
\begin{minipage}[t]{6cm}
\setlength{\unitlength}{0.07cm}
\begin{picture}(60,80)
\thinlines
\put(10,0)     {\framebox(60,80)}
\put(10,5)     {\line(1,0){60}}
\put(10,10)    {\line(1,0){60}}
\put(10,15)    {\line(1,0){60}}
\put(10,20)    {\line(1,0){60}}
\put(10,25)    {\line(1,0){60}}
\put(10,30)    {\line(1,0){60}}
\put(10,35)    {\line(1,0){60}}
\put(10,40)    {\line(1,0){60}}
\put(10,45)    {\line(1,0){60}}
\put(10,50)    {\line(1,0){60}}
\put(10,55)    {\line(1,0){60}}
\put(10,60)    {\line(1,0){60}}
\put(10,65)    {\line(1,0){60}}
\put(10,70)    {\line(1,0){60}}
\put(10,75)    {\line(1,0){60}}
\put(25,0)     {\line(0,1){80}} 
\put(40,0)     {\line(0,1){80}} 
\put(55,0)     {\line(0,1){80}} 
\put(10,0)     {\makebox(5,5){1}}
\put(10,75)    {\makebox(5,5){64}}

\linethickness{1mm}
\put(18,2.5)   {\line(0,1){15}}
\put(18,17.5)  {\line(1,0){14}}
\put(32,2.5)   {\line(0,1){15}}
\put(32,2.5)   {\line(1,0){30}}
\put(18,2.5)   {\line(0,1){15}}
\put(18,22.5)  {\line(0,1){35}}
\put(18,22.5)  {\line(1,0){14}}
\put(18,57.5)  {\line(1,0){14}}
\put(18,62.5)  {\line(0,1){15}}
\put(18,62.5)  {\line(1,0){14}}
\put(32,22.5)  {\line(0,1){15}}
\put(32,42.5)  {\line(0,1){15}}
\put(32,62.5)  {\line(0,1){15}}
\put(32,37.5)  {\line(1,0){30}}
\put(32,42.5)  {\line(1,0){30}}
\put(32,77.5)  {\line(1,0){30}}
\put(48,7.5)   {\line(1,0){14}}
\put(48,12.5)  {\line(1,0){14}}
\put(48,17.5)  {\line(1,0){14}}
\put(48,22.5)  {\line(1,0){14}}
\put(48,27.5)  {\line(1,0){14}}
\put(48,32.5)  {\line(1,0){14}}
\put(48,47.5)  {\line(1,0){14}}
\put(48,52.5)  {\line(1,0){14}}
\put(48,57.5)  {\line(1,0){14}}
\put(48,62.5)  {\line(1,0){14}}
\put(48,67.5)  {\line(1,0){14}}
\put(48,72.5)  {\line(1,0){14}}

\put(48,7.5)  {\line(0,1){5}}
\put(48,17.5) {\line(0,1){5}}
\put(48,27.5) {\line(0,1){5}}
\put(48,47.5) {\line(0,1){5}}
\put(48,57.5) {\line(0,1){5}}
\put(48,67.5) {\line(0,1){5}}
\put(62,2.5)  {\line(0,1){5}}
\put(62,12.5)  {\line(0,1){5}}
\put(62,22.5)  {\line(0,1){5}}
\put(62,32.5)  {\line(0,1){5}}
\put(62,42.5)  {\line(0,1){5}}
\put(62,52.5)  {\line(0,1){5}}
\put(62,62.5)  {\line(0,1){5}}
\put(62,72.5)  {\line(0,1){5}}
\end{picture}
\caption{The second iteration of method "P"}
\label{P-second iteration}
\end{minipage}
\hfill
\begin{minipage}[t]{6cm}
\setlength{\unitlength}{0.07cm}
\begin{picture}(60,80)
\thinlines
\put(10,0)     {\framebox(60,80)}
\put(10,5)     {\line(1,0){60}}
\put(10,10)    {\line(1,0){60}}
\put(10,15)    {\line(1,0){60}}
\put(10,20)    {\line(1,0){60}}
\put(10,25)    {\line(1,0){60}}
\put(10,30)    {\line(1,0){60}}
\put(10,35)    {\line(1,0){60}}
\put(10,40)    {\line(1,0){60}}
\put(10,45)    {\line(1,0){60}}
\put(10,50)    {\line(1,0){60}}
\put(10,55)    {\line(1,0){60}}
\put(10,60)    {\line(1,0){60}}
\put(10,65)    {\line(1,0){60}}
\put(10,70)    {\line(1,0){60}}
\put(10,75)    {\line(1,0){60}}
\put(25,0)     {\line(0,1){80}} 
\put(40,0)     {\line(0,1){80}} 
\put(55,0)     {\line(0,1){80}} 
\put(10,0)     {\makebox(5,5){1}}
\put(65,0)     {\makebox(5,5){64}}

\linethickness{1mm}
\put(18,2.5)   {\line(1,0){14}} 
\put(18,7.5)   {\line(1,0){14}} 
\put(18,12.5)  {\line(1,0){14}} 
\put(18,17.5)  {\line(1,0){14}} 
\put(18,22.5)  {\line(1,0){14}} 
\put(18,27.5)  {\line(1,0){14}} 
\put(18,32.5)  {\line(1,0){14}} 
\put(18,37.5)  {\line(1,0){14}} 
\put(18,42.5)  {\line(1,0){14}} 
\put(18,47.5)  {\line(1,0){14}} 
\put(18,52.5)  {\line(1,0){14}} 
\put(18,57.5)  {\line(1,0){14}} 
\put(18,77.5)  {\line(1,0){14}} 

\put(48,2.5)   {\line(1,0){14}} 
\put(48,7.5)   {\line(1,0){14}} 
\put(48,12.5)  {\line(1,0){14}} 
\put(48,17.5)  {\line(1,0){14}} 
\put(48,22.5)  {\line(1,0){14}} 
\put(48,27.5)  {\line(1,0){14}} 
\put(48,32.5)  {\line(1,0){14}} 
\put(48,37.5)  {\line(1,0){14}} 
\put(48,42.5)  {\line(1,0){14}} 
\put(48,47.5)  {\line(1,0){14}} 
\put(48,52.5)  {\line(1,0){14}} 
\put(48,57.5)  {\line(1,0){14}} 
\put(48,77.5)  {\line(1,0){14}} 
\put(32,62.5)  {\line(1,0){16}} 

\put(18,7.5)    {\line(0,1){5}} 
\put(18,17.5)   {\line(0,1){5}} 
\put(18,27.5)   {\line(0,1){5}} 
\put(18,37.5)   {\line(0,1){5}} 
\put(18,47.5)   {\line(0,1){5}} 
\put(18,57.5)   {\line(0,1){20}} 

\put(32,2.5)    {\line(0,1){5}} 
\put(32,12.5)   {\line(0,1){5}} 
\put(32,22.5)   {\line(0,1){5}} 
\put(32,32.5)   {\line(0,1){5}} 
\put(32,42.5)   {\line(0,1){5}} 
\put(32,52.5)   {\line(0,1){5}} 
\put(32,62.5)   {\line(0,1){15}} 

\put(62,7.5)    {\line(0,1){5}} 
\put(62,17.5)   {\line(0,1){5}} 
\put(62,27.5)   {\line(0,1){5}} 
\put(62,37.5)   {\line(0,1){5}} 
\put(62,47.5)   {\line(0,1){5}} 
\put(62,57.5)   {\line(0,1){20}} 

\put(48,2.5)    {\line(0,1){5}} 
\put(48,12.5)   {\line(0,1){5}} 
\put(48,22.5)   {\line(0,1){5}} 
\put(48,32.5)   {\line(0,1){5}} 
\put(48,42.5)   {\line(0,1){5}} 
\put(48,52.5)   {\line(0,1){5}} 
\put(48,62.5)   {\line(0,1){15}} 
\end{picture}
\caption{The second iteration of method "H"}
\label{H-second iteration}
\end{minipage}
\hfill
\end{figure}

\begin{definition}
 \rm Every $t\in [0,1]$ is uniquely determined by a sequence of nested closed intervals (that are generated by our successive partitioning),the lengths of which shrink to 0. With this sequence corresponds a unique sequence of nested closed squares, the diagonals of which shrink into a point, and which define a unique point in $[0,1]^2$, the image $S_n(t)$ of $t$. 
\end{definition}

\begin{theorem}
\label{big lemma}
\rm For every $n \in\mathbb N$ there exists a continuous cubes preserving function  $$f_n:[0,1]\stackrel{onto}{\longrightarrow} [0,1]^n.$$
\end{theorem}
\noindent
{\bf Proof.}
By the induction on $n\in\mathbb N$.
\newline For $n=1$ let $f_1$ be the identity map $I~:~[0,1]\rightarrow[0,1]$.
\newline For $n=2$ let $f_2$ be the Hilbert space-filling curve.
\newline Now if we have already defined $f_{n-1}:~[0,1]\rightarrow[0,1]^{n-1}$ continuous cubes preserving function, then let:
$$f_n=\left(\begin{array}{r} {f_{n-1}\circ \varphi}\\ \psi\end{array}\right)$$
where $\varphi,\psi$ are the coordinate functions of $S_n=\left(\begin{array}{r}\varphi\\\psi\end{array}\right)$ defined as in Lemma \ref{small lemma}.
\newline Then $f_n:~[0,1]\rightarrow[0,1]^n$ is continuous. Let us prove that $f_n$ is cubes preserving.\\

Let $\alpha\subseteq[0,1]$ be an element of $K^1_{n\cdot s}$  for some $s\in \mathbb N$. Then $S_n(\alpha) = \alpha^\prime\times\alpha^{\prime\prime}$  for some $\alpha^\prime\in K^1_{(n-1)\cdot s},~\alpha^{\prime\prime}\in K^1_s$.~~Or in other words~~$\varphi(\alpha) =\alpha^\prime,~\psi(\alpha)=\alpha^{\prime\prime}.$  It follows that $f_{n-1}\circ \varphi(\alpha)\subseteq\delta$  for some $\delta\in K^{n-1}_s$.  Consequently,
$$f_n(\alpha)=\left(\begin{array}{r} f_{n-1}\circ\varphi(\alpha)\\ \psi(\alpha)\end{array}\right)\subseteq\delta\times\alpha^{\prime\prime}\in K^n_s$$
So that the property (\ref{cubes preserving 1}) of cubes preserving function is proven.

Now let $\delta$  be an element of $K^n_s$  for some $s\in\mathbb N$.  Then $int(\delta)=int(\delta^\prime)\times int(\alpha)$  for some  $\delta^\prime\in K^{n-1}_s,~\alpha\in K^1_s.$
\newline We can present the function $f_n$ in matrix form as :
$$f_n=\left(\begin{array}{rr}f_{n-1} & 0 \\
                                   0 & I\end{array}\right)\cdot
      \left(\begin{array}{r}\varphi\\ \psi\end{array}\right)$$
or $f_n=F\cdot S_n$,  where $F=\left(\begin{array}{rr}f_{n-1} & 0 \\ 0 & I\end{array}\right)$  and $I$ is the identity map.
\newline Then $f_n^{-1}= S^{-1}_n\circ F^{-1}$. And we can see that
$$F^{-1}(int(\delta))\subseteq f^{-1}_{n-1}(int(\delta^\prime))\times I^{-1}(int(\alpha)).$$
By the induction hypothesis, we have $f^{-1}_{n-1} (int(\delta^\prime))\subseteq\alpha^\prime$  for some $\alpha^\prime\in K^1_{(n-1)\cdot s}$.  
\newline 
Recalling that $f_{n-1}$ is continuous and $int(\delta^\prime)$ is an open set, we can say that $f^{-1}_{n-1}(int(\delta^\prime))\subseteq int(\alpha^\prime)$. And obviously $I^{-1}(int(\alpha))=int(\alpha)$. So we can write $$F^{-1}(int(\delta))\subseteq int(\alpha^\prime)\times int(\alpha) =  int(\alpha^\prime\times\alpha)$$
where $\alpha^\prime\in K^1_{(n-1)\cdot s},~\alpha\in K^1_s$.

Further, by the property \ref{Sn property 2} of the function $S_n$ (see Lemma \ref{small lemma}) , it follows that\newline
$S^{-1}_n(int(\alpha^\prime\times\alpha))\subseteq\alpha^{\prime\prime}$ for some $\alpha^{\prime\prime}\in K^1_{n\cdot s}$. So finally we conclude
\newline
$f^{-1}_n(int(\delta))=S^{-1}_n(F^{-1}_n(int(\delta)))\subseteq S^{-1}_n(int(\alpha^\prime)\times int(\alpha))=S^{-1}_n(int(\alpha^\prime\times\alpha))\subseteq\alpha^{\prime\prime}\in K^1_{n\cdot s}$
The property (\ref{cubes preserving 2}) for the function $f_n$ is proven, and hereby the theorem \ref {big lemma} is proven as well.\\
$\hspace*{\fill}\Box$

\begin{theorem}
\label{corollary}
\rm For every $n \in\mathbb N$ there exists a continuous cubes preserving function  $$f_n:[0,1]\stackrel{onto}{\longrightarrow} [0,1]^n.$$
with the property: \newline\noindent
if $[a,b],[b,c]\in K^1_{ns}$ for some $s\in\mathbb N$, then $f_n([a,b]),f_n([b,c])\in K^n_s$ have $(n-1)$-dimensional face of side $\frac{1}{2^s}$ in common, and $f_n(b)$ is a vertex of the face. 
\end{theorem}
\noindent
{\bf Proof.}
The proof is the same as the proof of the Theorem \ref{big lemma} and the property is just a property of the special functions $f_n$ constructed in the proof of the Theorem \ref{big lemma}.\\
$\hspace*{\fill}\Box$

\section{Proof of the Theorem \ref{theorem 1}.}
Let us define the function $p$. Our plan is first to define $p$ on some closed set $B_0\subseteq [0,1]^m$~such that $p(B_0) = [0,1]^n$. Then we will extend the definition of $p$ over all cube $[0,1]^m$ with the property that $p$ on $[0,1]^m $ is a $C^{<\frac{m}{n}}$-function and every partial derivative of order $<\frac{m}{n}$~vanishes on $B_0$.
Finally joining appropriate pairs of points of $B_0$ by straight line segments, each mapping by $p$ on a singular point in $[0,1]^n$, we get a connected set $E$ such that $B_0\subset E\subset [0,1]^m$~ and ~$p(E)=[0,1]^n.$~ Showing that every partial derivative of order $<\frac{m}{n}$~ vanishes on $E$, we will finish the proof of the Theorem \ref{theorem 1}.\\

For the construction of function $p$, we enumerate a set of subcubes of $[0,1]^n$~ and a set of subcubes of $[0,1]^m.$ In case of the cube $[0,1]^n$, a set of subcubes is the $K^n=\bigcup\limits_{s\in\mathbb N}K^n_s$~ with the enumeration $\langle K^n,\prec_{f_n}\rangle$ (see (\ref{linear order})), where $f_n$~is as in Theorem \ref{corollary} for the number $n\in\mathbb N$. So that for every $s\in\mathbb N, ~K^n_s=\{Q^n_{1,i_1,i_2,\dots,i_j,\dots,i_s};~1\leqslant i_j\leqslant 2^n\}$ where $Q^n_{1,i_1,i_2,\dots,i_s}~\prec_{f_n}~Q^n_{1,i_1^\prime,i_2^\prime,\dots,i_s^\prime}$~if~$\nexists j_0~(1\leqslant j_0\leqslant s)~$such that$~\forall j<j_0~~i_j\geqslant i_j^\prime$~and~$i_{j_0}> i_{j_0}^\prime$,~and $\forall i_s~(1\leqslant i_s\leqslant 2^n)~~ Q^n_{1,i_1,i_2,\dots,i_{s-1},i_s}\subset Q^n_{1,i_1,i_2,\dots,i_{s-1}}$.\\

\noindent
In case $[0,1]^m$ a set of subcubes, which we designate by $\tilde K^m=\bigcup\limits_{s\in \mathbb N}\tilde K^m_s$,~ where $\forall s\in\mathbb N,~\tilde K^m_s=\{\tilde Q^m_{1,i_1,i_2,\dots,i_j,\dots,i_s}~;~1\leqslant i_j\leqslant i_s\}$, is constructed similarly to the standard set $K^m=\bigcup\limits_{s\in\mathbb N}K^m_s$ as follows:
\begin{itemize}
\item $\tilde Q^m_1=[0,1]^m$~ is the cube of side $S_1=1=\frac{1}{2^0}$
\item $\{\tilde Q^m_{1,i_1}~;~i_1=1,2,\dots,2^m\}$~ is a set of $m$-dimensional subcubes of $\tilde Q^m_1$~ of side $S_2=\frac{1}{2^1}(1-\frac{1}{2^2})=\frac{1}{2}S_1(1-\frac{1}{4})=\frac{3}{8}<\frac{1}{2}S_1$,~ such that in every corner of $\tilde Q^m_1$~ there locates the only one subcube $\tilde Q^m_{1,i_1}$~for some $i_2\leqslant 2^m.$
\item Now if the cube $\tilde Q^m_{1,i_1,i_2,\dots,i_{s-1}}$~ have already been defined for some $s\in\mathbb N$,~then we define a set of $m$-dimensional subcubes $\tilde Q^m_{1,i_1,i_2,\dots,i_{s-1},i_s}~\subset~\tilde Q^m_{1,i_1,i_2,\dots,i_{s-1}}~;~i_s=1,2,\dots,2^m$,~ each of side $S_s=\frac{1}{2^{s-1}}\prod\limits_{j=1}^s (1-\frac{1}{j^2})=\frac{1}{2}S_{s-1}(1-\frac{1}{s^2}) < \frac{1}{2}S_{s-1}$,~such that in every corner of $\tilde Q^m_{1,i_1,i_2,\dots,i_{s-1}}$~ there locates the only one subcube $Q^m_{1,i_1,i_2,\dots,i_{s-1},i_s}$ for some $i_s\leqslant 2^m$.
\end{itemize} 

Because of similarity of the construction of the sets $K^m$ and $\tilde K^m$, the standard enumeration of the set $K^m$~:~$\langle K^m,\prec_{f_m}\rangle$,~where $f_m$~ is as  in the Theorem \ref{corollary} for the number $m\in\mathbb N$,~~ can be easily transfered on the set $\tilde K^m$. Let us designate that enumeration by $\langle \tilde K^m,\prec_{f_m}\rangle$.\\
The figures \ref{order on K^2},\ref{order on tilde K^2} show for case $m=2$ the enumeration $\langle K^m,\prec_{f_m}\rangle$ on the set $K^m$~~ and the corresponding enumeration $\langle \tilde K^m,\prec_{f_m}\rangle$ on the set $\tilde K^m$ for the second iteration. Recall that for $m=2~~f_m=f_H$~-the Hilbert space-filling function.\\

\begin{figure}[ht]
\begin{minipage}[t]{6.5cm}
\setlength{\unitlength}{0.08cm}
\begin{picture}(60,60)
\put(0,0) {\framebox(60,60){$Q^2_1$}}
\put(0,0)   {\framebox(30,30){$Q^2_{11}$}}
\put(0,30)  {\framebox(30,30){$Q^2_{12}$}}
\put(30,0)  {\framebox(30,30){$Q^2_{14}$}}
\put(30,30) {\framebox(30,30){$Q^2_{13}$}}
\put(0,0)    {\framebox(15,15){$Q^2_{111}$}}
\put(0,15)   {\framebox(15,15){$Q^2_{114}$}}
\put(15,0)   {\framebox(15,15){$Q^2_{112}$}}
\put(15,15)  {\framebox(15,15){$Q^2_{113}$}}

\put(30,30)  {\framebox(15,15){$Q^2_{131}$}}
\put(30,45)  {\framebox(15,15){$Q^2_{132}$}}
\put(45,30)  {\framebox(15,15){$Q^2_{134}$}}
\put(45,45)  {\framebox(15,15){$Q^2_{133}$}}

\put(30,0)    {\framebox(15,15){$Q^2_{143}$}}
\put(30,15)   {\framebox(15,15){$Q^2_{142}$}}
\put(45,0)    {\framebox(15,15){$Q^2_{144}$}}
\put(45,15)   {\framebox(15,15){$Q^2_{141}$}}

\put(0,30)    {\framebox(15,15){$Q^2_{121}$}}
\put(0,45)    {\framebox(15,15){$Q^2_{122}$}}
\put(15,30)   {\framebox(15,15){$Q^2_{124}$}}
\put(15,45)   {\framebox(15,15){$Q^2_{123}$}}
\end{picture}
\caption{$\langle K^2,\prec_{f_2}\rangle$~on $K^2$}
\label{order on K^2}
\end{minipage}
\hfill
\begin{minipage}[t]{6.5cm}
\setlength{\unitlength}{0.08cm}
\begin{picture}(60,60)
\put(0,0)   {\framebox(60,60){$\tilde Q^2_1$}}
\put(0,0)   {\framebox(25,25){$\tilde Q^2_{11}$}}
\put(0,35) {\framebox(25,25){$\tilde Q^2_{12}$}}
\put(35,0)  {\framebox(25,25){$\tilde Q^2_{14}$}}
\put(35,35) {\framebox(25,25){$\tilde Q^2_{13}$}}
\put(0,0)   {\framebox(10,10){$\tilde Q^2_{111}$}}
\put(0,15)  {\framebox(10,10){$\tilde Q^2_{114}$}}
\put(15,0)  {\framebox(10,10){$\tilde Q^2_{112}$}}
\put(15,15) {\framebox(10,10){$\tilde Q^2_{113}$}}

\put(35,35) {\framebox(10,10){$\tilde Q^2_{131}$}}
\put(35,50) {\framebox(10,10){$\tilde Q^2_{132}$}}
\put(50,35) {\framebox(10,10){$\tilde Q^2_{134}$}}
\put(50,50) {\framebox(10,10){$\tilde Q^2_{133}$}}

\put(35,0)  {\framebox(10,10){$\tilde Q^2_{143}$}}
\put(35,15) {\framebox(10,10){$\tilde Q^2_{142}$}}
\put(50,0)  {\framebox(10,10){$\tilde Q^2_{144}$}}
\put(50,15) {\framebox(10,10){$\tilde Q^2_{141}$}}

\put(0,35)  {\framebox(10,10){$\tilde Q^2_{121}$}}
\put(0,50)  {\framebox(10,10){$\tilde Q^2_{122}$}}
\put(15,35) {\framebox(10,10){$\tilde Q^2_{124}$}}
\put(15,50) {\framebox(10,10){$\tilde Q^2_{123}$}}
\end{picture}
\caption{$\langle\tilde K^2,\prec_{f_2}\rangle$~on $\tilde K^2$}
\label{order on tilde K^2}
\end{minipage}
\end{figure}

Now we turn to the definition of the function $p:[0,1]^m\rightarrow[0,1]^n$. We establish some conditions the function $p$ must satisfy. By induction: $p(\tilde Q^m_1)= Q^n_1$, and if 
\begin{align*}
&p(\tilde\delta)=\delta~~~~~\text{for}~~~ \tilde\delta\in\tilde K^m_{n(s-1)},~\delta\in K^n_{m(s-1)}
\end{align*}
then sets $\tilde\delta^*=\{\tilde\gamma\in\tilde K^m_{ns};~\tilde\gamma\subseteq\tilde\delta\}
,~~\delta^*=\{\gamma\in K^n_{ms};~\gamma\subseteq\delta\}$~ each have $2^{nm}$ elements, and both are linear ordered by $\prec_{f_m},~\prec_{f_n}$~ respectively. So that we have an order preserving bijection between the sets $\tilde\delta^*$~and~$\delta^*$. For every $\tilde\gamma\in\tilde\delta^*$~ let us write 
\begin{align*}
&p(\tilde\gamma)=\gamma\in\delta^*~~\text{such that~~if}~~\tilde\gamma,\tilde\gamma^\prime\in\tilde\delta^*,~~~\tilde\gamma\prec_{f_m}\tilde\gamma^\prime~~\Longrightarrow~~p(\tilde\gamma)\prec_{f_m}p(\tilde\gamma^\prime).
\end{align*}

Now let us define the set $B_0\subseteq \tilde Q_1^m$ and $p\upharpoonright B_0$~ such that $p(B_0)= [0,1]^n$. 
Let $B_0$  be a set of points $x\in \tilde Q_1^m$,~ each of which is uniquely determined by a sequence of nested closed $m$-dimensional cubes (that are generated by our successive partitioning), diagonals of which shrink into the point.

 Let $x$ be an element of $B_0$~ and $\{P_i~:~i\in\mathbb N\}$ be the sequence of the $m$-dimensional cubes $\tilde Q^m_{1i_1\cdots i_{j_n}}$ that determine the point~$x$, then $x= \bigcap\limits_{i\in\mathbb N} P_i$~~and also $\bigcap\limits_{i\in\mathbb N}p(P_i)\in [0,1]^n$. Let us define $p(x)=\bigcap\limits_{i\in\mathbb N} p(P_i)$. It is not difficult to see that $p:B_0\stackrel{onto}{\rightarrow}[0,1]^n$.

 We have completed definition of the set $B_0\subseteq [0,1]^m$~and the function $p:B_0\rightarrow [0,1]^n$.   Let $p_r:B_0\rightarrow[0,1];~r=1,\dots,n$ be the coordinate function for $p:B_0\rightarrow [0,1]^n$. Then $\forall t\in B_0$
$$p(t)=(p_1(t),p_2(t),\dots,p_n(t))\in [0,1]^n.$$
To construct a $C^{<\frac{m}{n}}$ extension of $p$ over all $[0,1]^m$, we need to find a $C^{<\frac{m}{n}}$ extension of $p_r$ over $[0,1]^m$ for every $r~~~(1\leqslant r\leqslant n).$ Let us fix such a $p_r$ for some $r~~~(1\leqslant r\leqslant n).$

To extend the function $p_r$, we will need to  introduce an intermediate function $p_{x^\prime,x^{\prime\prime}}$. 
Let $x^\prime=(x_1^\prime\cdots x_j^\prime \cdots x^\prime_m), 
~x^{\prime\prime} =(x_1^{\prime\prime}\cdots x_j^{\prime\prime}\cdots x_m^{\prime\prime})\in \mathbb R^m$~be such that: 
\begin{align}
\label{intermediate function}  
&\exists j^\prime\leqslant m~:~\forall j\not=j^\prime~~~~x^\prime_j =x^{\prime\prime}_j \text{~~and~~}x^\prime_{j^\prime}\leqslant x^{\prime\prime}_{j^\prime}.
\end{align}
For such $x^\prime,~x^{\prime\prime}$~ let us define $L_{x^\prime,x^{\prime\prime}}=\{x\in\mathbb R^m~:~x_j=x^\prime_j~\text{if}~j\not = j^\prime ~\text{and}~ x^\prime_j\leqslant x_j\leqslant x_j^{\prime\prime}\}$.
And $\forall x\in L_{x^\prime,x^{\prime\prime}}$~we define 
\begin{align}
\label{Munkres function}
&p_{x^\prime,x^{\prime\prime}} (x) = (p_r(x^{\prime\prime}) - p_r(x^\prime))\cdot g\left(\frac{x_i-x_i^\prime}{x_i^{\prime\prime}-x_i^\prime}\right)+ p_r(x^\prime)
\end{align}
where (following \cite{Munkres},p.6) $g:\mathbb R\rightarrow[0,1]$~ is a smooth map such that 
\begin{align*}
&g\upharpoonright(-\infty,0] = 0,\\
&g\upharpoonright[1,\infty) = 1,\\
&g^\prime(t)> 0 ~~\text{for}~~ 0<t<1.
\end{align*}

 Therefore $p_{x^\prime,x{\prime\prime}}:~L_{x^\prime,x^{\prime\prime}}\rightarrow [p_r(x^\prime),p_r(x^{\prime\prime})]\subseteq \mathbb R$. Obviously, $p_{x^\prime,x^{\prime\prime}}$~makes sense only if $p_r(x^\prime),~p_r(x^{\prime\prime})$ have already been defined. Then $p_{x^\prime,x^{\prime\prime}}$~ is a smooth, strictly monotone increasing bijection if $p_r(x^\prime)< p_r(x^{\prime\prime})$, or is a constant if $p_r(x^\prime)= p_r(x^{\prime\prime})$
, such that 
\begin{align*}
 &p_{x^\prime,x^{\prime\prime}}(x^\prime)= p_r(x^\prime),\\
 &p_{x^\prime,x^{\prime\prime}}(x^{\prime\prime}) = p_r(x^{\prime\prime}),\\
 &(p_{x^\prime,x^{\prime\prime}})^{(l)}_{x_i}(x^\prime)=  (p_{x^\prime,x^{\prime\prime}})^{(l)}_{x_i}(x^{\prime\prime})= 0~~\forall l\in\mathbb N,\\
 &(p_{x^\prime,x^{\prime\prime}})^\prime_{x_i} (x) > 0~~\forall x\in L_{x^\prime,x^{\prime\prime}}\setminus (x^\prime\cup x^{\prime\prime})~~~~\text{if}~~p_r(x^\prime)<p_r(x^{\prime\prime})\\
 &(p_{x^\prime,x^{\prime\prime}})^\prime_{x_i} (x) = 0~~\forall x\in L_{x^\prime,x^{\prime\prime}}\setminus (x^\prime\cup x^{\prime\prime})~~~~\text{if}~~p_r(x^\prime)=p_r(x^{\prime\prime}).
\end{align*}

 We will define $p_r$~on $[0,1]\setminus B_0$ sequentially. As for all $x\in B_0~~p_r(x)$ have already been defined, we begin by definition $p_r(x)$~ for all $x$ that lie on some edge of some cube $\tilde Q^m_{1i_1\dots i_{jn}},~j\in\mathbb N$.~ Let us designate the set of such $x$ as $B_1\subseteq [0,1]^m$.~ For any $x\in B_1\setminus B_0$~ there exits $L_{x^\prime,x^{\prime\prime}}\ni x$,~ lying on the same as $x$ edge of some cube $\tilde Q^m_{1i_1\dots i_{jn}}$, and $L_{x^\prime,x^{\prime\prime}}\cap B_0 = \{x^\prime,x^{\prime\prime}\}$~ because, by the construction, $B_0$ is a closed set containing all vertexes of all cubes in $\tilde K^m_{jn}$.   For $x\in B_1\setminus B_0$~ we define $p_r(x) = p_{x^\prime,x^{\prime\prime}}(x)$.

  Now by induction, let us define a set $B_l\subseteq [0,1]^m$~ as a set of all points $x$ lying on some $l$-dimensional face of some cube $\tilde Q^m_{1i_1\cdots i_{jn}}\in \tilde K^m_{jn}$, and $x\not\in B_{l-1}$.~~ Suppose that we have defined $p_r(x),~x\in\bigcup\limits_{i=0}^{l-1}B_i$. If $x\in B_l\setminus B_0$, then there exists $\tilde Q(x)\in\tilde K^m_{jn}$~ such that $x=(x_1\cdots x_m)\in \tilde Q(x),~x\not \in \tilde Q(x)_{i_{jn+1}\dots i_{(j+1)n}}$.

If there exists $L_{x^\prime,x^{\prime\prime}}\ni x$,~where $x^\prime,x^{\prime\prime}$~satisfy \eqref{intermediate function} and  
\begin{align*}
&{ \left.\begin{array}{r}
x^\prime\in B_{l-1}\cap \tilde Q(x)_{i^\prime_{jn+1}\dots i^\prime_{(j+1)n}}\\
x^{\prime\prime}\in B_{l-1}\cap \tilde Q(x)_{i^{\prime\prime}_{jn+1}\dots i^{\prime\prime}_{(j+1)n}}
  \end{array}
 \right\}~~\text{for~some}~~i^\prime_{jn+1}\dots i^\prime_{(j+1)n}\not = i^
{\prime\prime}_{jn+1}\dots i^{\prime\prime}_{(j+1)n}}\\
&int (L_{x^\prime,x^{\prime\prime}})\subseteq \tilde Q(x)\setminus \bigcup\limits_{i_{jn+1},\dots,i_{(j+1)n}=1}^{2^m} \tilde Q(x)_{i_{jn+1}\dots i_{(j+1)n}}
\end{align*}
then $p_r(x) \stackrel{def}{=} p_{x^\prime,x^{\prime\prime}}(x).$

  Let us designate all points $x$ of the set $B_l$ for which such $L_{x^\prime,x^{\prime\prime}}$~~ exists  as $B^\prime_l$. One can notice that $B^\prime_1=B_1\setminus B_0$.\\

  We can suppose that $p_r(x)$ is already defined for all $x\in B^\prime_l$. And now
for every $x\in B^{\prime \prime}_l = B_l\setminus (B^\prime_l \cup B_0)$ let us define $p_r(x) = p_{x^\prime,x^{\prime\prime}}$,~~ where $x^\prime,~x^{\prime\prime}\in \tilde Q(x)\cap B_l^\prime $  satisfy (\ref{intermediate function}) with least possible $j^\prime$~~required by the (\ref{intermediate function}),~~ such that  $L_{x^\prime,x^{\prime\prime}}\ni x$~~  and~~ $int(L_{x^\prime,x^{\prime\prime}})\subseteq \tilde Q(x)\setminus \bigcup\limits_{i_{jn+1},\dots,i_{(j+1)n}=1}^{2^m} \tilde Q(x)_{i_{jn+1}\dots i_{(j+1)n}}$.\\

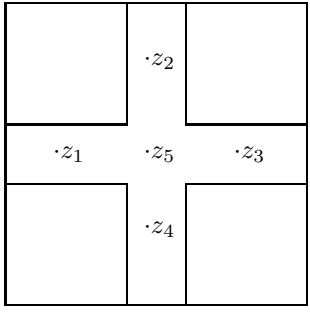
\setlength{\unitlength}{0.08cm}
\begin{picture}(65,65)
\put(0,0)  {\framebox(50,50)}
\put(0,0)   {\framebox(20,20)}
\put(0,30)  {\framebox(20,20)}
\put(30,0)  {\framebox(20,20)}
\put(30,30) {\framebox(20,20)}

\put(25,10) {\makebox(1,5){$\cdot z_4$}}
\put(25,38) {\makebox(1,5){$\cdot z_2$}}
\put(10,22.5) {\makebox(1,5){$\cdot z_1$}}
\put(40,22.5) {\makebox(1,5){$\cdot z_3$}}
\put(25,22.5) {\makebox(1,5){$\cdot z_5$}}
\put(70,27) {\makebox{This picture shows an example for $m=2$}}
\put(70,22.5)  {\makebox{$z_1,z_2,z_3,z_4\in B_2^\prime,~~z_5\in B_2\setminus B_2^\prime$ .}}
\end{picture}\\

As $p\upharpoonright B_0$ was defined earlier, the definition of $p_r:[0,1]^m\rightarrow [0,1]$ has been completed.  The definition of the function $p$ is well defined as it follows from the fact that the family of  sets $\{B_l;~l\leqslant m\}$,~~ is disjoint and for any $x\in \bigcup\limits_{l=1}^m B_l \setminus B_0$~~ there exists a unique $L_{x^\prime,x^{\prime\prime}}$.\\

 For every $x\in [0,1]^m\setminus B_0$,~~ let us designate the unique set $L_{x^\prime,x^{\prime\prime}}\ni x$~~ as $L(x)$ ~~ and $L^\prime (x)\stackrel{def}{=} x^\prime, ~~L^{\prime\prime}(x)\stackrel{def}{=} x^{\prime\prime}$,~~and if $x\in B_0$,~~then $L(x)=L^\prime(x)=L^{\prime\prime}(x)=x$.\\

To show that $p_r\in C^{<\frac{m}{n}}$ and every partial derivative of $p_r$ of order $<\frac{m}{n}$ vanishes on $B_0$ we need to prove the following three lemmas.
  
\begin{lemma}\label{inner lemma 1}
\rm 
$$\lim\limits_{j(s)\to\infty}\frac{\diam(p_r(\tilde Q^m_{1i_1\dots i_{j(s)n}}))}{\left(S_{(j(s)+1)n-1} - 2S_{(j(s)+1)n} \right)^k}=0.$$
\end{lemma}

\noindent
{\bf Proof.}
\begin{align*}
&\lim\limits_{j(s)\to\infty}\frac{\diam(p_r(\tilde Q^m_{1i_1\dots i_{j(s)n}}))}{\left(S_{(j(s)+1)n-1} - 2S_{(j(s)+1)n} \right)^k}~\leqslant\\
&\lim\limits_{j(s)\to\infty}\frac{\sqrt{n}\left(\frac{1}{2}\right)^{j(s)m-1}}{\left(\frac{1} {2^{(j(s)+1)n-2}}\cdot\frac{1}{((j(s)+1)n)^2}\cdot\prod\limits_{i=2}^{(j(s)+1)n-1}(1-\frac{1}{i^2}) \right)^k}~=\nonumber
\end{align*}
\begin{align}
&\sqrt{n}\cdot 2^{(n-2)k+1}\cdot\lim\limits_{j(s)\to\infty}\frac{((j(s)+1)n)^{2k}}{2^{j(s)n(\frac{m}{n}-k)}}\cdot \lim\limits_{j(s)\to\infty}\frac{1}{\left(\prod\limits_{i=2}^{(j(s)+1)n-1}(1-\frac{1}{i^2})\right)^k}~=~0\label{line 4}
\end{align}

which follows from the fact that the first limit in (\ref{line 4}) is equivalent to $\lim\limits_{n\to\infty}\frac{n^\alpha}{a^n}=0$,~for $\alpha>0,~a>1$;~ and $\lim\limits_{s\to\infty}\prod\limits_{i=2}^s(1 - \frac{1}{i^2})\geqslant b > 0$~~for some $b\in \real^+$, because the series $\sum \frac{1}{i^2}$~~converges \cite{Knopp}.\\
Lemma \ref{inner lemma 1} is proven.\\

\begin{lemma}\label{inner lemma 2}\rm
If $\{x^s;~s\in\mathbb N\}\subseteq [0,1]^m\setminus B_0$,~such that $\lim\limits_{s\to\infty}|L(x^s)|=0$,~then $\forall k~~(0\leqslant k<\frac{m}{n})$
$$\lim\limits_{s\to\infty}\frac{p_r(L^{\prime\prime}(x^s)) - p_r(L^\prime(x^s))}{|L(x^s)|^k}=0.$$
\end{lemma}

\noindent
{\bf Proof.}
\newline\noindent
If $\tilde Q(x^s)= \tilde Q^m_{1i_1\dots i_j(s)n}$,~~ then 
\begin{align*}
&\lim\limits_{s\to\infty}\frac{|p_r(L^{\prime\prime}(x^s))-p_r(L^\prime(x^s))|}{|L(x^s)|^k}~\leqslant\\
&\lim\limits_{j(s)\to\infty}\frac{\diam(p_r(\tilde Q^m_{1i_1\dots i_{j(s)n}}))}{\left(S_{(j(s)+1)n-1} - 2S_{(j(s)+1)n} \right)^k}~=~0~~\text{by Lemma \ref{inner lemma 1}.}
\end{align*}
Lemma \ref{inner lemma 2} is proven.\\

\begin{lemma}\label{inner lemma 3}\rm
Let $x^0=(x_1^0,x_2^0,\dots,x_i^0,\dots,x_m^0)\in B_0$~ and~ $k<\frac{m}{n}$,~~then \newline\noindent
$\lim\limits_{x\to x_0}\frac{p_r(x)-p_r(x^0)}{|x_i-x_i^0|^k}=0$~ for $x=(x_1^0,x_2^0,\dots,x_i,\dots,x_m^0).$
\end{lemma}

\noindent
{\bf Proof.}
\newline\noindent
To prove that, suffice it to show that the limit is equal to 0 for some arbitrary chosen sequence  $\{x^j=(x_1^0,x_2^0,\cdots,x_i^j,\cdots,x_m^0)\in B^\prime_{i_0}\setminus B_0\}$,~~ also without loss of generality we can suppose that $x_i^j > x_i^{j+1}$. 
Let us define $\forall j\in \mathbb N~~j(s)$~ as the lowest number such that $x^j,x^0\in \tilde Q^m_{1,j_2,\dots,i_{j(s)n}}$~ for some $\tilde Q^m_{1,j_2,\dots,i_{j(s)n}}\in \tilde K^m_{j(s)n}$, then 
$$\lim\limits_{x\to x_0}\frac{p_r(x)-p_r(x^0)}{|x_i-x_i^0|^k} \leqslant 
\lim\limits_{j(s)\to\infty}\frac{\diam(p_r(\tilde Q^m_{1i_1\dots i_{j(s)n}}))}{\left(S_{(j(s)+1)n-1} - 2S_{(j(s)+1)n} \right)^k}=0~~ \text{by Lemma \ref{inner lemma 1}}.$$
Lemma \ref{inner lemma 3} is proven.\\

Note: if $f\in C^k$,~then it is not difficult to see that $f^{(k^\prime)}\equiv 0~~~\forall k^\prime\in \real^+\setminus \mathbb Z,~k^\prime < k.$\\

We will use the induction by $i~~~(1\leqslant i\leqslant m).$ \\
For $i_0=1,~p_r\in C^\infty$~ on $B_1\setminus B_0$~ and if $x\in B_1\setminus B_0,~k>0,$~ then
$$D_kp_r(x)=\frac{p_r(x^{\prime\prime})-p_r(x^\prime)}{(x_i^{\prime\prime}-x^\prime_i)^k}\cdot D_kg\left(\frac{x_i-x_i^\prime}{x_i^{\prime\prime}-x_i^\prime}\right) $$\newline
where $x^\prime=L^\prime(x),~x^{\prime\prime}=L^{\prime\prime}(x),~x_i$~ is unique coordinate that is not a constant on $L(x)$. If  $x\in B_1\cap B_0$,~then~$D_kp_r(x)=0$ for every $k<\frac{m}{n}$. For the proof see {\bf Case 2} for $i_0=1$.\\

Now let $i_0\in\mathbb N~~(i_0\leqslant m)$~~ be such that $p_r\upharpoonright\bigcup\limits_{i<i_0}B_i$~~ is a $C^{<\frac{m}{n}}$-function and all its partial derivatives of order less than $\frac{m}{n}$~~ vanish on~~ $\bigcup\limits_{i<i_0}B_i\cap B_0$.

\noindent
Let us consider some $k~~~(0\leqslant k < \frac{m}{n})$,and a set $T(k)=\{0,1,\cdots,[k],k, k-1,\cdots,k-[k]\}$.\\

\noindent
Let $D_tp_r(x)$~~ be  short for $\frac{\partial^{t_1+\cdots+t_m}}{\partial x_1^{t_1}\cdots\partial x_m^{t_m}} p_r(x_1^\prime\cdots x^\prime_m).$ \\

\noindent
{\bf Case 1}. \newline
If $x\in \bigcup\limits_{i\leqslant i_0}B^\prime_i\setminus B_0$,~~ then considering derivatives with respect to every variable that is not constant in a neighborhood of  $x\in \bigcup\limits_{i\leqslant i_0}B^\prime_i\setminus B_0$, we write the following general formula
\begin{align}
\label{derivative}
&D_{k}p_r(x) = \frac{D_{k-t}p_r(x^{\prime\prime})- D_{k-t}p_r(x^\prime)} {(x_i^{\prime\prime}-x_i^\prime)^t} D_tg\left(\frac{x_i-x_i^\prime}{x_i^{\prime\prime}-x_i^\prime}\right) + D_{k}p_r(x^\prime)
\end{align}
for some $t\in T(k)$.

\noindent
where $x^\prime=L^\prime(x),~~x^{\prime\prime}=L^{\prime\prime}(x),~~x_i$~~is the unique coordinate that is not a constant on $L(x).$\\

\noindent
$D_{k}p_r$~~exists at point $x$~~ and is continuous. It follows from the choice of $i_0$ and from $g$~~ being a $C^\infty$-function.\\

\noindent
We also can notice that if $t\not = 0$,  then $D_{k}p_r(x^\prime)= 0$;~ and if~$t$~is not integer,~then $D_tg=0.$\\

\noindent
{\bf Case 2}.\newline  
Let $x^0 \in \bigcup\limits_{i\leqslant i_0}B^\prime_i\cap B_0$. Without loss of generality, we can suppose that $p_r\upharpoonright \bigcup\limits_{i\leqslant i_0}B^\prime_i$~~ is a $C^{k^\prime}$-function,~~~  where 
\begin{equation}
 k^\prime = \left\{
\begin{array}
{l@{\quad \text{if} \quad}l}
[k] & k \text{~~is~not~ integer}\\
k-1 & k \text{~~is~integer},
\end{array}\right.
\label{assumption on k}
\end{equation} 
\noindent
and all partial derivatives of order $\leqslant k^\prime$~~vanish on $B_0\cap \bigcup\limits_{i\leqslant i_0}B^\prime_i$.\\

\noindent
Now we need to prove that 
\begin{equation}
\label{label2}
D_{k}(x^0) = \lim\limits_{x\to x^0}\frac{|D_{k^\prime}p_r(x) - D_{k^\prime}p_r(x^0)|}{|x_i-x_i^0|^\lambda}=0
\end{equation}

\noindent
where $x=(x_1^0,x_2^0,\cdots,x_i,\cdots,x_m^0)$~~ for  some $i$ depending on the $D_{k}$,\newline 
and 
$ \lambda = \left\{
\begin{array}
{l@{\quad \text{if} \quad}l}
k-[k] & k \text{~~is~not~ integer}\\
1 & k \text{~~is~integer}.
\end{array}\right.$ 

To prove that, suffice it to show that the limit is equal to 0 for some arbitrary chosen sequence  $\{x^j=(x_1^0,x_2^0,\cdots,x_i^j,\cdots,x_m^0)\in B^\prime_{i_0}\setminus B_0\}$,~~ also without loss of generality we can suppose that $x_i^j > x_i^{j+1}$.  Then for every point $x^j~~~(j\in\mathbb N)$~~ by the definition of $L(x^j)$~~ and that $x^0\in B_0$, it follows that $L(x^j)$~~ lies on the straight line $(x_1^0,x_2^0,\cdots,x_i,\cdots,x_m^0)$~~ and 
\begin{equation}
L^\prime(x^j),L^{\prime\prime}(x^j)\in B_0\cap\bigcup\limits_{i\leqslant i_0}B_i^\prime.
\label{condition 1}
\end{equation}

\noindent
Substituting the sequence $\{x^j;~j\in\mathbb N\}$ in~~(\ref{label2}),  we write:
\begin{align}
&\lim\limits_{x^j\to x^0}\frac{|D_{k^\prime}p_r(x^j) - D_{k^\prime}p_r(x^0)|}{|x_i^j-x_i^0|^\lambda}~~\leqslant\nonumber\\
&\lim\limits_{x^j\to x^0}\frac{|D_{k^\prime}p_r(x^j) - D_{k^\prime}p_r(L^\prime(x^j))|}{|x_i^j-L^\prime(x^j)_i|^\lambda}~~=\label{label8}\\
&\lim\limits_{x^j\to x_0}\frac{|D_{k^\prime-t}p_r(L^{\prime\prime}(x^j)) - D_{k^\prime-t}p_r(L^\prime(x^j))|}{|L^{\prime\prime}(x^j)_i-L^\prime(x^j)_i|^t}\cdot \frac{D_tg\left(\frac{x^j_i - L^\prime(x^j_i)}{L^{\prime\prime}(x^j)_i - L^\prime(x^j)_i}\right)}{|x^j_i - L^\prime(x^j)_i|^\lambda}~~=\label{label9}\\
&\lim\limits_{x^j\to x_0}\frac{|D_{k^\prime-t}p_r(L^{\prime\prime}(x^j)) - D_{k^\prime-t}p_r(L^\prime(x^j))|}{|L^{\prime\prime}(x^j)_i-L^\prime(x^j)_i|^t}\nonumber\\
&\times\frac{D_tg\left(\frac{x^j_i - L^\prime(x^j_i)}{L^{\prime\prime}(x^j)_i - L^\prime(x^j)_i}\right) - D_tg\left(\frac{L^\prime(x^j_i)- L^\prime(x^j_i)}{L^{\prime\prime}(x^j)_i-L^\prime(x^j)_i}\right)}{\left|\frac{|x^j_i-L^\prime(x^j)_i|}{L^{\prime\prime}(x^j)_i-L^\prime(x^j)_i}-\frac{L^\prime(x^j_i)-L^\prime(x^j_i)}{L^{\prime\prime}(x^j)_i-L^\prime(x^j)_i}\right|^\lambda}\cdot\frac{1}{(L^{\prime\prime}(x^j)_i-L^\prime(x^j)_i)^\lambda}~~=\label{label5}\\
&\lim\limits_{x^j\to x_0}\frac{|D_{k^\prime-t}p_r(L^{\prime\prime}(x^j)) - D_{k^\prime-t}p_r(L^\prime(x^j))|}{|L^{\prime\prime}(x^j)_i-L^\prime(x^j)_i|^{t+\lambda}}\cdot
\frac{D_tg(y)-D_tg(0)}{|y-0|^\lambda} ~~=0\label{label6}
\end{align} 

where
\begin{itemize}
\item in (\ref{label2}):  $D_{k^\prime}p_r(x^0) = 0$~~ by our assumption on $k^\prime$ (see (\ref{assumption on k}));
\item in (\ref{label8}): $D_{k^\prime}p_r(L^\prime(x^j))=0$~~if $k^\prime >0$, because $L^\prime(x^j) \in B_0\cap \bigcup\limits_{i\leqslant i_0}B_i^\prime$~~(when $k^\prime=0$~~ the limit is equal to zero by Lemma \ref{inner lemma 3});
\item in (\ref{label9}): see (\ref{derivative});
\item in (\ref{label5}): $D_tg(0)=0~~~~\forall t\in T$;
\item substitution in (\ref{label6}):  $y=\frac{x^j_i-L^\prime(x^j_i)}{L^{\prime\prime}(x^j)_i-L^\prime(x^j)_i}$~~;
\item in case that there exists $j_0\in\mathbb N$~~such that $\forall j\geqslant j_0~~~L^\prime(x^j)=x_0$,~~limit of the first factor in (\ref{label6}) is bounded  and limit of the second factor is equal to 0. Otherwise the second factor in (\ref{label6}) is bounded and the limit of the first factor is equal to 0 by the Lemma \ref{inner lemma 2} if $k^\prime-t =0$~~~or by (\ref{condition 1})  and (\ref{assumption on k})~~if $k^\prime-t \not = 0$.
\end{itemize}

\noindent
{\bf Case 3}.\newline  
For $x^0\in \bigcup\limits_{i\leqslant i_0}(B_i\setminus B_i^\prime)$~~ the proof is similar to {\bf Case 1}, assuming that $p_r\upharpoonright (\bigcup\limits_{i<i_0}B_i \bigcup B_{i_0}^\prime)$ is a $C^{k}$-function.\\

Now it follows from the cases 1,2,3:
\newline\noindent
if $p_r\upharpoonright\bigcup\limits_{i<i_0}B_i$~ is $C^{<\frac{m}{n}}$~ function and all its partial derivatives of order $<\frac{m}{n}$ vanish on $\bigcup\limits_{i<i_0}B_i\cap B_0$,~ then $p_r\upharpoonright\bigcup\limits_{i\leqslant i_0}B_i$~ is $C^{<\frac{m}{n}}$ function and  all  its partial derivatives of order $<\frac{m}{n}$~ vanish on $\bigcup\limits_{i\leqslant i_0}B_i\cap B_0$. So that by induction, we have that $p_r\upharpoonright[0,1]^m$~ is $C^{<\frac{m}{n}}$~ function and all its partial derivatives of order $<\frac{m}{n}$ vanish on $B_0$. \\

\noindent
{\bf The arc.}\newline
If~ $x(1,i_1,\dots,i_s)=x(1,i_1,\dots,i_{s+1})=b$~is a common vertex of cubes 
\begin{align*}
&Q^m_{1,i_1,i_2,\dots,i_s}=f_m([a,b]),~~Q_{1,i_1,i_2,\dots,i_{s+1}}^m=f_m([a,b])\in K^m_s
\end{align*}
for some $s\in \mathbb N$, and $[a,b],[b,c]\in K^1_{ns}$, then the corresponding to it vertexes $$\tilde x(1,i_1,\dots,i_s),~~\tilde x(1,i_1,\dots,i_{s+1})~~\text{of cubes}~~ 
\tilde Q^n_{1,i_1,i_2,\dots,i_s},~~\tilde Q^m_{1,i_1,i_2,\dots,i_{s+1}}\in \tilde K^m_s$$
respectively are, by the construction of $\tilde K^m$,~ elements of $B_0$ satisfying  (\ref{intermediate function}),~ and, by the definition of function $p$~on~$B_0$,
$$p(\tilde x(1,i_1,\dots,i_s))=p(\tilde x(1,i_1,\dots,i_{s+1})).$$
\noindent 
So that there exists the $L_{\tilde x(1,i_1,\dots,i_s),\tilde x(1,i_1,\dots,i_{s+1})}$~ and from (\ref{Munkres function})\newline 
$p\upharpoonright L_{\tilde x(1,i_1,\dots,i_s),\tilde x(1,i_1,\dots,i_{s+1})}$~ is a constant. And more than that, because of the fact that 
$\tilde x(1,i_1,\dots,i_s)$,~$\tilde x(1,i_1,\dots,i_{s+1})\in B_0$~ all partial derivatives of $p$ of order $<\frac{m}{n}$~vanish at those points and, as follows from (\ref{derivative}), they vanish on \newline
$L_{\tilde x(1,i_1,\dots,i_s),\tilde x(1,i_1,\dots,i_{s+1})}$.\\

Let us define $E=B_0\cup L_{\tilde x(1,i_1,\dots,i_j,\dots,i_s),\tilde x(1,i_1,\dots,i_j,\dots,i_{s+1})},~(1\leqslant i_j\leqslant 2^m,~s\in\mathbb N).$~~ Then $E\subseteq[0,1]^m,~p(E)=[0,1]^n,$~ every partial derivative of $p$ of order $<\frac{m}{n}$~ vanishes on $E$,~and the connectivity of $E$ can be shown by an argument similar to that of Whitney \cite{Whitney}.\\
\hspace*{\fill}$\Box$ \\

\section{Proof of the Theorem \ref{theorem 2}.}
Using the result of the Theorem \ref{theorem 1} we can take a map $p^\prime:[0,1]^{m-r}\rightarrow[0,1]^{n-r}$~ contained in $C^k$~ for all real $k<\frac{m-r}{n-r}$, and a connected subset $E^\prime$~ such that every partial derivative of $p^\prime$~ of order $<\frac{m-r}{n-r}$~ vanishes on $E^\prime$~ and $p^\prime(E^\prime)=[0,1]^{n-r}$. Then $p=(I,p^\prime)$~ and $E=E^\prime\times[0,1]^r$,~ where $I:[0,1]^r\rightarrow[0,1]^r$~is the identity map, satisfy the conditions of the Theorem \ref{theorem 2}. \\
\hspace*{\fill}$\Box$ \\

\end{document}